# TEXTUAL ANALYSIS OF ANCIENT INDIAN MATHEMATICS

Satyanad Kichenassamy

**Abstract.** Recent analyses of Brahmagupta's discourse on the cyclic quadrilateral, and of Baudhāyana's approximate quadrature of the circle, have shown that it is useful to submit mathematical texts to a form of literary analysis. Several passages considered as obscure or objectionable may be explained in this way, by taking into account the elements of exposition and derivation of the results that the author has given, as well as his conceptual background. This approach aims at helping the reader set aside his preconceptions about what a mathematical text is supposed to be. In this paper, guidelines for further application of this method are outlined, with illustrations taken from our previous papers.

Indian mathematical discourse, as we find it in the earliest primary sources, has often been criticized for apparently leaving definitions, derivations, or even some of the hypotheses to be supplied by the teacher when the student enquires about them. Now, in the case of propositions XII.21-32 of Brahmagupta's *Brāhmasphuṭasiddhānta*[i] (*BSS*), and propositions I.58-62 of Baudhāyana's *Śulvasūtra*[ii] (*BŚlS*) we have recently shown[iii] that these texts contain more than meets the eye: they do form consistent discourses that contain precise statements as well as elements of derivations. We have suggested[iv] that the key to the understanding of these texts is to perform a literary analysis. Indeed, the text should not be taken as a sequence of bits of information, a list of results, but as a structured discourse, that aims at conveying an argument. To understand this argument, we suggest one should not only consider *what* the author says, but also *how* he says it, and investigate *why* he said it so and not otherwise.

This paper is organized as follows. In Sect. 1, the difficulties posed by these two texts are recalled. In Sect. 2, we show that literary analysis is useful here because it helps setting aside preconceptions about the text. In Sect. 3, we outline a few examples; we analyze a neologism used by Brahmagupta (§ 3.1), the expression of perpendicularity in *BSS* XII.30-31 (§ 3.2), the quadrature rule *BŚlS* I.59 (§ 3.3), and an apparent redundancy in *BSS* XII.24 (§ 3.4). The reader is referred to the papers already quoted for further examples and background information, including a discussion of earlier attempts at accounting for these passages. We conclude (Sect. 4) with a few common-sense guidelines[v] that may help in the application of this approach to further texts.

---

[i] 628 AD, see (Colebrooke, 1817; Dvivedin, 1902; Sharma, 1966).
[ii] First half of the first millennium BC, see (Thibaut, 1875; Sen & Bag, 1983).
[iii] See (Kichenassamy, 2006; Kichenassamy, 2010; Kichenassamy, 2012a).
[iv] See (Kichenassamy, 2012b).

## 1. Difficulties in construing Indian mathematical texts.

The two texts we consider have been found difficult to construe for different reasons. Brahmagupta's discourse on the cyclic quadrilateral, *BSS* XII.21-32, was criticized for apparently not specifying the cyclicity condition nor, when needed, the condition of orthogonality of the diagonals. Recall that this passage is a sequence of statements consisting of one or several verses each, in the Sanskrit metre *āryā*.[vi] It describes, in a very specific terminology, without figures or numerical examples, general properties of lines in cyclic quadrilaterals, with occasional hints on the logical relation between them. It includes, apparently for the first time, celebrated results such as the expressions of the area and diagonals in terms of the sides, or the result known as "Brahmagupta's theorem" in modern mathematics.[vii] We suggested that the very structure of this discourse contains the apparently missing information.[viii] We do not dwell here on the tools that Brahmagupta had at his disposal; we merely point out that he does not use the notions of angle and parallel, not even that of right angle, but that he knows how to drop or erect perpendiculars; the roots for these constructions are to be found in the *Śulva*s.

Baudhāyana's *BŚlS* Prop. I.58-62[ix] pose a different problem. They include a complex and unusual rule for the approximate quadrature of the circle (I.59), found in no other text, in India or elsewhere. It appears at face value that it is based on the manipulation of fractions, as we find it described in later Indian mathematical texts; but it has not been possible to account for it in this manner, despite systematic attempts. The issues may be described as follows.[x] Calling *r* the radius of the circle, and *s* the side of the square with the same area, I.58 expresses an approximate circulature rule: if the side is given, the radius is taken to be

$$r = \frac{s}{2} + \frac{E}{3}$$

where $E = (\sqrt{2} - 1)\frac{s}{2}$ is the excess of the half-diagonal over the half-side. Thus, the application of the circulature rule depends on the knowledge of an approximation of $\sqrt{2}$. I.59 and I.60 are quadrature rules; the symmetry of the formulation of I.58 and I.59-60 suggests that the quadrature rules were obtained by inversion of the circulature rule. I.59 appears to be equivalent to the modern statement (in which *d* stands for the diameter):

---

[v] Many hermeneutic traditions are implicitly making use of similar guidelines. Our suggestion is that such methods may be adapted so as to become relevant to the understanding of scientific texts.
[vi] In this paper, each such statement will be called a proposition.
[vii] See § 3.2 below for the statement of this result.
[viii] See (Kichenassamy, 2010; Kichenassamy, 2012a).
[ix] In the numbering of (Thibaut, 1875); they correspond respectively to 2.8-12 in (Sen & Bag, 1983).
[x] The following remarks are due to Cantor, Müller and Drenckhahn, see (Kichenassamy, 2006).

$$s = d\left\{1 - \frac{1}{8 \times 29}\left[28 + \frac{1}{6} - \frac{1}{8} \times \frac{1}{6}\right]\right\}.$$

I.60 gives the cruder rule $s = d\{1 - 2/15\}$. Finally, I.61-62 gives an approximation of $\sqrt{2}$:

$$1 + \frac{1}{3} + \frac{1}{4} \times \frac{1}{3} - \frac{1}{34} \times \left(\frac{1}{4} \times \frac{1}{3}\right).$$

If we assume that the reduction of fractions to the same denominator was used by Baudhāyana, this approximation takes the form 577/408. Substituting into the expression for the radius given in (I.58), we obtain

$$r = \frac{1393}{1224} \times \frac{s}{2}.$$

This yields $d/s = 1393/1224$. Inverting, we obtain $s/d = 1224/1393$. Now, this value is very close to the one given in I.59, since

$$\frac{1224}{1393} = \left\{1 - \frac{1}{8 \times 29}\left[28 + \frac{1}{6} - \frac{1}{8} \times \frac{1}{6}\right]\right\} - \frac{41}{8 \times 29 \times 6 \times 8 \times 1393}.$$

Thus, the approximation I.61-62 does not quite yield I.59, but a very close value. This "verification" is not a derivation of I.59 because the decomposition of 1224/1393 is not unique:

$$\frac{1224}{1393} = \left\{1 - \frac{33}{8 \times 34}\right\} + \frac{1}{8 \times 34 \times 1393}.$$

This formula involves division by 34, just as the approximation of $\sqrt{2}$. Thus, the assumption that Baudhāyana manipulated fractions does not account for his text. We need to answer the following questions: what is the origin of the division into 8 parts, followed by division into 29 parts? Why did Baudhāyana state that one should remove 28 parts? In other words, why did he not state the apparently simpler rule

$$s = d\left\{1 - \frac{1}{8} + \frac{1}{8 \times 34}\right\}?$$

## 2. The relevance of literary analysis to the study of scientific texts.

Literary analysis of poems or similar works is a natural endeavor, because they use the resources of ordinary language to express very specific, but perhaps unfamiliar ideas or feelings. The words are taken from the common stock, but may be used in novel ways, so that some thought may be needed before their exact meaning is unraveled. Also, a poem may have several dimensions, may be read at several levels. Thus, part of the analysis of a poem consists in collecting "contextual information" that was obvious to the author, but not necessarily to us. But there is another issue: if the text does not correspond to our expectations, we may misunderstand its import. It is common experience that difficulties in construing poems, even in one's mother tongue, may be impaired by one's *a priori* conception of poetic discourse.[xi] Literary analysis has developed in part in order to provide safeguards against the intrusion of one's own preconceptions.

Scientific texts share these characteristics of poems. They use ordinary words in a technical sense, but remain within the framework of the common structures of language. They aim, in the case of original works, at conveying new results. Nevertheless, they are not always immediately or fully intelligible. Questions of priority or the like may occasionally intrude and cause confusion or obscurity, but that is not the difficulty we have in mind. We all know that modern mathematicians are prone to call "trivial" a statement for which any graduate student ought to be able to supply proof. In other words, it is tacitly assumed that, in any given branch of mathematics, there is a perspective from which certain results become obvious, and that graduate instruction is meant to communicate this point of view. Put otherwise, graduate instruction instills certain preconceptions that make a class of results "plain." Since all mathematicians do not share the same preconceptions, they may not, unless they know each other rather well, understand each other without some amount of analysis. We need to perform a similar analysis when studying ancient mathematicians, because we do not know beforehand which results they considered too obvious for inclusion in an advanced treatise such as Brahmagupta's.

### 3.  Extracts from analyses of mathematical texts.

We now illustrate these issues by four examples taken from our analyses of *BSS* XII.21-32 and *BŚlS* I.58-62.

### *3.1. Brahamgupta's "triquadrilateral"*

Proposition XII.21 of *BSS* states his famous formula for the area of a cyclic quadrilateral:

$$Area = \sqrt{(s-a)(s-b)(s-c)(s-d)},$$

where $s = (a+b+c+d)/2$ is the half-sum of the sides *a*, *b*, *c* and *d*. It is often stated that Brahamgupta did not specify that his quadrilaterals were required to be inscribed in a circle. Now, Brahmagupta does not refer to a quadrilateral, but uses the neologism *tricaturbhuja*, that may be rendered by "triquadrilateral." The above criticism of Brahmagupta stems from the assumption that "triquadrilateral" means "triangle and (unrelated) quadrilateral." Let us examine whether this assumption is consistent with the text. This word is used in *BSS* in two propositions: XII.21 and XII.27. It does not seem to be used in any other Sanskrit work. Let us infer what figure the triquadrilateral must be, on the sole basis of the two propositions in *BSS* where the term appears.

  **XII.21**  *sthūlaphalaṃ tricaturbhujabāhupratibāhuyogadalaghātaḥ*
        *bhujayogārdhacatuṣṭayabhujonaghātāt padaṃ sūkṣmam*

    *Crude is the value (of the area) given by*
    *The product of the half-sums of opposite sides of a triquadrilateral;*
    *From a group of four half-sums of the sides,*
    *Having subtracted (each) side (in turn), the foot*[xii] *of the product is the refined (value).*

---

[xi] Richards' experiments illustrate this point; see (Richards, 1929).

**XII.27** *tribhujasya vadho bhujayordviguṇitalamboddhṛto hṛdayarajjuḥ*
*sā dviguṇā tricaturbhujakoṇaspṛgvṛttaviṣkambhaḥ*

*Of the trilateral, the product of both sides,*
*Divided by twice the perpendicular, is the heart-cord;*
*The same, doubled, is, of the triquadrilateral,*
*The diameter of the circle that touches its corners.*

The "triquadrilateral" must have four sides, because XII.21 mentions four. It does not refer to an unrelated triangle and quadrilateral, because XII.27 would then be redundant: it would state the expression of the circumradius of triangle twice. In addition, XII.21 would need to be modified to apply to a triangle: since a triangle has only three sides, say $a$, $b$ and $c$, it would be necessary to replace the fourth term $s - d$ by $s$. There is no textual evidence for the statement, occasionally made, that Brahmagupta considered a triangle as a quadrilateral of vanishing fourth side. XII.27 states that there is a circle touching the corners of the triquadrilateral, and that its diameter is twice the circumradius of the triangle considered in its first hemistich. Therefore, the triquadrilateral contains a distinguished triangle. In order to find the circumradius of this triangle, we need to know its sides, and one of them is a diagonal of the quadrilateral. However, the diagonals of the cyclic quadrilateral have not been determined at this point of the text (they are given in XII.28). Therefore, the only reasonable conclusion is that *the triquadrilateral is the figure obtained by completing a triangle by the addition of a fourth vertex on its circumcircle*.

As was pointed out by P.-S. Filliozat, our analysis is supported by the standard procedures for the formation of compounds described by Sanskrit grammar: *tricaturbhuja* is a *bahuvrīhi* (or possessive compound) of which the first member is itself a *dvandva* (or copulative compound): it should be decomposed as [*tri-catur*]-*bhuja*.[xiii] It therefore refers to a single object, possessing "sides" that are "three and four" in number.[xiv]

Returning to XII.21, we recall that Brahmagupta gives not only the exact formula for the area of a cyclic quadrilateral, but also a very gross one: the product of the half-sums of opposite sides. What is the point of giving such a formula? Now, unlike the exact formula, the gross formula makes reference to opposite sides: it depends on the order in which the sides are taken. If we view the cyclic quadrilateral as obtained from a circle by removing four segments of a circle, just as, in the *Śulva*s, a square differs from its circumcircle by four segments (called *pradhi*s), we see that the area of a quadrilateral determined by four consecutive chords in a circle cannot depend on the order of its sides, since a given chord determines the areas of the two segments it bounds, without reference to any other chord. We therefore suggested that Brahmagupta was stressing the independence of the area from the order of the sides, criticizing implicitly the gross formula. This independence is useful in the derivation of the exact formula.

---

[xii] The "foot" (*pada*) refers to the (square) root (*mūla*). Indeed, the roots of a tree spread at its foot.
[xiii] "We propose to analyse *tricaturbhuja* as a *bahuvrīhi* made of *bhuja* qualified by the two terms of the *dvandva tricatur*" (Filliozat, 2010b).
[xiv] There are three sides if we view it as a triangle to which a fourth vertex is added, four sides if we view it as a quadrilateral. As in the famous metaphor of the elephant and the blind, a single object may be apprehended differently depending on one's focus of interest. This is consistent with the meaning of the *dvandva* compound: "Besides its most common duty of expressing coordination, the dvandva is also available if "and" connects persons or things standing in mutual relation with one another." (Speijer, 1886, §207; this reference also gives examples from the classics).

## 3.2. Perpendicularity without angles

*BSS* XII.30-31 is the source of the result called "Brahmagupta's theorem" in modern textbooks. Its modern statement is as follows: in a cyclic quadrilateral with diagonals that cut each other at right angles, the perpendicular to one side, through the intersection of the diagonals, bisects the opposite side. This formulation is due to Chasles (Chasles, 1837). Brahmagupta states an equivalent result that will not be recalled here.[xv] It seems to be found in no other work before 1837, apart from commentaries on *BSS,* in India or elsewhere. Now, Brahmagupta does not use the notion of angle, not even that of right angle. However, he routinely uses the notion of perpendicular dropped from a vertex of a triangle or quadrilateral to a side, as well as those of "segment" (*āvādhā*) and "portion" (*khaṇḍa*). Consider for instance the triangle bounded by the two adjacent sides *a* and *b*, and diagonal γ of a cyclic quadrilateral; we take γ as base of the triangle (Fig. 1). The foot *H* of the perpendicular of this triangle, dropped to γ, divides it into two parts α and β, called segments. They are the projections of the two adjacent sides on the base. Similarly, the point *J,* where the diagonals of the quadrilateral meet, divides γ into two parts, called portions. In XII.30-31, Brahmagupta expresses that, in the two triangles based on the diagonals, the segments are also portions: *tadāvādhe pṛthagūrdhvādharakhaṇḍe* (XII.30). This means that the perpendicular coincides with part of the second diagonal. Therefore, in modern terms, the diagonals are orthogonal. We conclude that Brahmagupta expresses the condition of perpendicularity of the diagonals by the identity of segments and portions.[xvi] This is a convenient formulation because Brahmagupta has no term to express directly that two lines "cut each other at right angles."

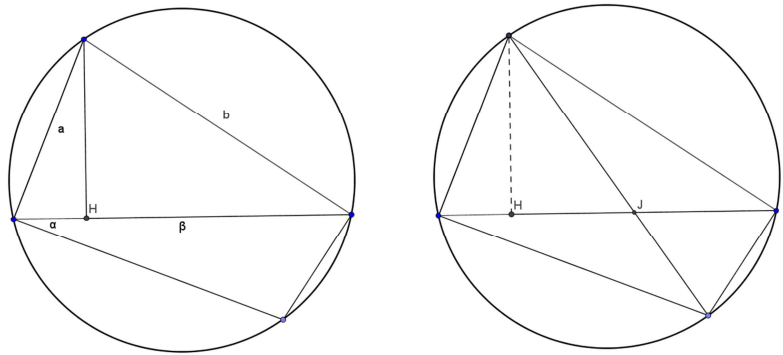

**Figure 1.** Definition of segments and portions. Left: one diagonal is taken as base of a triangle. The foot *H* of the perpendicular dropped to the base determines two segments *α* and *β* on it. Right: the intersection *J* of the diagonals determines two portions on each of them. Segments and portions coincide if *H=J*, so that diagonals "cut each other at right angles."

## 3.3. "Scale-calculus" in Baudhāyana's Śulva.

Consider *BŚlS* I.59:

**I.59.** *maṇḍalaṃ caturaśraṃ cikīrṣanviṣkambhamaṣṭau bhāgānkṛtvā bhāgamekona-triṃśadhā vibhajyāṣṭāviṃśatibhāgānuddharet |*

---

[xv] See (Kichenassamy, 2012a).
[xvi] See (Kichenassamy, 2012a; Kichenassamy, 2012b).

> *bhāgasya ca ṣaṣṭhamaṣṭamabhāgonam*
>
> (Someone) wishing to turn a circle into a square, having made out of the diameter eight parts, (and) divided a part twenty-nine-fold, should remove twenty-eight (of these twenty-nine) parts;
>
> and (moreover) the sixth (part) of the (one) part (left) less the eighth part (of the sixth).

As we recalled in Sect. 1, earlier results indicate that the quadrature rule I.59 represents the result of an inversion of the circulature rule I.58, and that a calculus of fractions cannot account for the way this inversion was carried out. We have therefore suggested that Baudhāyana used a different mathematical notion, which we called *scale-calculus*.[xvii] We outline some of its features, with application to I.59.

Here as in all of *BŚlS*, the manipulation of lengths that may be divided into an arbitrary number of parts plays a central role. These lengths are geometrically represented by marked cords. As one sees from the text, the correspondence between the diameter of a circle and the equivalent square is expressed by the decomposition of one length into a number of equal parts, say *p*, such that *q* of these parts give another length. There is no reference here to an absolute unit of length: all we have is a relation between two lengths. This point is essential in the enlargement of figures that is performed routinely in the *Śulva*s: it is achieved by taking a new unit, defined geometrically and not arithmetically. In other words, the new unit is not necessarily defined as an integral number of parts of the old unit, but may be constructed. Another example is given by I.58: the circulature is given as a correspondence between the excess $E = (\sqrt{2} - 1)s/2$ of the half-diagonal over the half-side, and the excess $e = r - s/2$ of the radius over the half-side. Prop. I.58 expresses that $e = E/3$. Here, the new unit $E$ is not given as a rational multiple of the side. Thus, the operations in the text should be described in terms of a correspondence (*p, q*) with a scalable unit of length. Such a correspondence between two lengths is readily inverted: if we divide the second length into *q* parts, *p* of these yield the first length. Here, *p* and *q* play symmetric roles, unlike the numerator and denominator of a fraction.

Unfortunately, I.58 is not phrased as a correspondence between side and diameter, but between $E$ and $e$, both of which involve the side of the square. If this side is known, this readily yields the radius, hence the diameter. But if the radius is given and not the half-side, it is not clear how to determine the latter. To determine a (*p, q*) correspondence between side and diameter amounts to finding some new unit of length subdivided into a certain number of parts, such that the side contains *p* parts, and the diameter *q*. For instance, take as unit the half-diagonal, divided into twelve parts. It turns out that I.59 may be accounted for as the result of the improvement of the crude (upper and lower) approximations $1\frac{1}{2}$ and $1\frac{1}{3}$ for $\sqrt{2}$. These two approximations respectively correspond to a half-side of 8 or 9 parts. A natural improvement is to take $8\frac{1}{2}$. In other words, if the half-diagonal is divided into 24 parts, the half-side contains 17. Applying I.58 yields that if the half-side contains 17 parts, the radius is longer by one-third of seven parts. Dividing each part into three further parts, we obtain from this a $(51 = 3 \times 17, 7)$ correspondence between the half-side and excess *e*. Therefore, the diameter contains $58 = 51 + 7$ parts. Conversely, if the diameter is divided into 58 parts, the side is obtained by removing 7 parts. Since 58 is close to 56, which is a multiple of 7, this (58,

---

[xvii] We use "calculus" in the sense of "a particular method or system of calculation or reasoning." No connection with infinitesimal calculus is implied.

7) correspondence is very close to (56, 7), that is, to (8, 1). To make use of this, it is natural to further subdivide each part into four: the (58, 7) correspondence becomes $(4 \times 58, 4 \times 7) = (8 \times 29, 28)$. We have therefore accounted for the division, first into 8, and then into 29 parts, as well as for the removal of 28 parts. Observe how Baudhāyana, by the formulation of his rule, has given hints that constrain possible derivations. For further details, including a derivation of the second part of I.59, and of I.60-62, see Sect. 5 of (Kichenassamy, 2006).

### *3.4. Did Brahmagupta indulge in redundancy?*

Any apparent redundancy in a text written by a careful author must make us wonder whether the impression that some content is duplicated is not the result of a faulty interpretation. This point is particularly important in the case of Indian Mathematics, where authors strive for extreme brevity. As P.-S. Filliozat points out:[xviii] "If the interpreter encounters in a statement a form which he can prove superfluous, he must either conclude that the formulation is faulty or carry on with his effort to understand until he finds the information that may legitimately be derived from the statement brought into question."

An example of an apparent redundancy is provided by *BSS* XII.24, where the result now known as "Pythagoras' theorem" appears to be stated three times with simple permutations of the terms.

> **XII.24.** *karṇakṛteḥ koṭikṛtiṃ viśodhya mūlam bhujaḥ bhujasya kṛtim*
> *prohya padaṃ koṭiḥ koṭibāhukṛtiyutipadaṃ karṇaḥ*
>
> *From the diagonal squared the upright squared subtracted,*
> *The root (of the difference) is the side. The square of the side*
> *Removed, the foot is the upright.*
> *The root of the sum of upright and arm is the diagonal.*

This appears to state, in a half-oblong with side (or arm) *a,* upright *b* and diagonal *c,* the relations

$$\sqrt{c^2 - b^2} = a; \sqrt{c^2 - a^2} = b; \sqrt{b^2 + a^2} = c.$$

They have already been used: XII.22 and 23 follow from them. In addition, the first two have exactly the same content, since the sides of an oblong play symmetric roles. It is odd that someone as keen on brevity as Brahmagupta should indulge in redundancy such as this. However, this interpretation tacitly assumes that the three statements in this proposition are three theorems about the same right triangle (or rather, half-oblong). This assumption is not necessary. A more satisfactory interpretation is to take them to refer to different half-oblongs. Consider a triangle *ABC* with sides *a=AB* and *b=BC*, of which the base $\gamma = AC$ is a diameter of the circumcircle: $\gamma = 2r$, where *r* is its circumradius (Fig. 2). Our assumption states that the midpoint *O* of the base is the circumcenter: $OA = OB = OC = r$. We do not know at this stage that *ABC* is a half-oblong. This figure contains several other triangles that are half-oblongs; they all have the perpendicular *h=BH* of the triangle as upright. Two of them, namely *ABH* and *BHC*, have as bases the segments $\alpha = AH$ and $\beta = HC$ determined by the perpendicular. The third, *OBH*, has its diagonal equal to the radius; its base is $OH = r - \alpha$.

---

[xviii] (Filliozat, 2010a) p. 27; our translation.

Note that $\beta = 2r - \alpha$, and $\gamma = \alpha + \beta$. Applying the first part of XII.24 to *ABH* and *BHC*, we obtain:

$$\alpha = \sqrt{a^2 - h^2}; \ \beta = \sqrt{b^2 - h^2}.$$

Applying the second part of XII.24 to *OBH*, we obtain

$$h^2 = r^2 - (r - \alpha)^2 = \alpha(2r - \alpha) = \alpha\beta.$$

Therefore, $\gamma^2 = (\alpha + \beta)^2 = \alpha^2 + \beta^2 + 2\alpha\beta = (a^2 - h^2) + (b^2 - h^2) + 2h^2 = a^2 + b^2$.

We have recovered the third statement in XII.24. The conclusion is that a triangle of which the base is a diameter of the circumcircle is necessary a half-oblong. Thus, the first two statements of XII.24 give the steps of the derivation of a result stated in the third. XII.24 with this interpretation is useful in the derivation of part of XII.26.

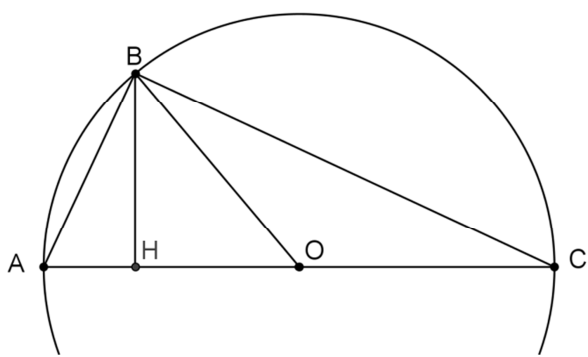

**Figure 2 :** Triangle having as base *AC* a diameter of the circumcircle of *ABC*. Its perpendicular is *BH*. It contains the half-oblongs *ABH*, *BHC*, and *OBH*.

### 4. Guidelines for the analysis of mathematical texts.

The preceding examples illustrate possible errors created by preconceptions; we conclude this paper with a few suggestions that may help detect them:

1. Unusual expressions should alert the reader: neologisms, apparent redundancies and the like, call for analysis.
2. The order of statements should be studied. If a statement seems to disrupt the logical order, one should investigate why.
3. The wording of statements should be studied. Why was this word used in preference of another seemingly equivalent one?
4. Operations used in any analysis should be limited to those allowed by the author's background, and preferably to operations attested in the text itself.
5. The mathematical consequences of assumptions about the *modus operandi* of the author should be followed to their end.
6. The structure of the text should be considered as a reflection of the author's view of the subject.
7. Each text should be taken on its own terms. Even a commentary, let alone a later independent text, may be a clarification, refutation or continuation of the text, or an opinion on it—it is no substitute for the text.

The above suggestions are meant to help pay attention to the small-scale structure of the text and to its consistency in the large. They may be summarized into the rule "Never assume." It is usually the consistency requirement that is the strongest tool to detect inadequate assumptions. One should also be aware that Indian texts such as the *Brāhmasphuṭasiddhānta* do not correspond to the norm of mathematical discourse widely disseminated through elementary modern teaching. They do not seek to establish that all results ultimately rely on a small number of common notions, or to present to a passive reader arguments to which he has no choice but to submit. On the contrary, the reader is expected to be active. It now appears that such texts are an invitation to think with the author, to rediscover his results with his help. Thus, the text is composed so as to stress the logical connection of one statement with the next, enabling the reader to follow the author, step by step. The common background notions are omitted as too elementary: it would be in bad taste, in an assembly of scholars, to dwell on such matters. Unfortunately, the partial breach in the continuity of Indian Mathematics, evidenced by the difficulty of later authors in understanding Brahmagupta, requires us to infer these elementary notions as well. That is another reason why a close analysis of these texts is required.

**Contact Details:**

**Satyanad Kichenassamy**
Laboratoire de Mathématiques
Université de Reims Champagne-Ardenne
Moulin de la Housse, B.P. 1039
F-51687 Reims Cedex 2, France.

*E-mail* : satyanad.kichenassamy@univ-reims.fr